\documentclass[11pt]{article}
\usepackage{a4,bm,epsfig,booktabs,pstricks}
\usepackage{setspace}
\usepackage{ifthen,array}
\newcommand{\ams}{\usepackage{amsfonts,amssymb,amsmath}}

\ams
\allowdisplaybreaks[3]
\newlength{\textwidthorig}
\newlength{\oddsidemarginorig}
\newlength{\textheightorig}
\newlength{\topmarginorig}
\def\seitenlaengenabsolut#1 #2 #3 #4 {\setlength{\textwidth}{#1}
                                      \setlength{\oddsidemargin}{#2}
                                      \setlength{\textheight}{#3}
                                      \setlength{\topmargin}{#4}}
\def\seitenlaengenrelzustandard#1 #2 #3 #4 {\setlength{\textwidth}{\textwidthorig+#1}
                                            \setlength{\oddsidemargin}{\oddsidemarginorig+#2}
                                            \setlength{\textheight}{\textheightorig+#3}
                                            \setlength{\topmargin}{\topmarginorig+#4}}
\def\seitenlaengenrelzuvorher#1 #2 #3 #4 {\addtolength{\textwidth}{#1}
                                          \addtolength{\oddsidemargin}{#2}
                                          \addtolength{\textheight}{#3}
                                          \addtolength{\topmargin}{#4}}
\newcommand{\standardseite}{\seitenlaengenrelzuvorher2.2cm -0.8cm 1.8cm -1.5cm }   %
\standardseite
\newcommand{\leerezeile}{\vspace{2ex}}
\newlength{\laengespatium}

\newcommand{\nach}{\longrightarrow}      

\newcommand{\auf}{\longmapsto}           
\newcommand{\txtauf}[1]{\auf}            
\newcommand{\impliz}{\Longrightarrow}    
\newcommand{\aequ}{\Longleftrightarrow}  
 
\newcommand{\invimpliz}{\Longleftarrow}  
\newcommand{\gegen}{\rightarrow}         
\newcommand{\iso}{\cong}                 
\newcommand{\ident}{\equiv}              
\newcommand{\teilmenge}{\subseteq}       
\newcommand{\obermenge}{\supseteq}       
\newcommand{\aeqrel}{\sim}               

\newcommand{\leeremenge}{\varnothing}    
\newcommand{\kreuz}{\times}              

\newcommand{\einschr}[1]{{}\arrowvert_{#1}}      
\newcommand{\dirsum}{\oplus}           

\newcommand{\betraganpass}[1]%
           {\left| #1 \right|}           
\newcommand{\bigbetrag}[1]%
           {\bigl|{#1}\bigr|}            
\newcommand{\betrag}[1]%
           {|{#1}|}                      
\newcommand{\betragnichtanpass}[1]%
           {\mid #1 \mid}                
\newcommand{\norm}[1]%
           {{}{\parallel}#1{\parallel}{}}      
\newcommand{\erww}[1]%
           {\langle #1 \rangle}          
\newcommand{\skalprod}[2]%
           {\langle #1,#2 \rangle}       
\newcommand{\supnorm}[1]{{\norm{#1}_\infty}}        
\newcommand{\quer}{\overline}            
\newcommand{\NULL}{\mathbf{0}}                         
\newcommand{\EINS}{{\boldsymbol{1}}}                   
\newcommand{\field}[1]{\mathbb{#1}}                    
\newcommand{\C}{{\field{C}}}                           
\newcommand{\R}{{\field{R}}}                           
\newcommand{\boundfkt}{{\cal B}}                       
\newcommand{\rnkl}[2]{\raisebox{-0.4ex}{$#1$}%
\raisebox{-0.12ex}{{\large$\setminus$}}\,#2}   
\newcommand{\agb}{{\overline{{\cal A}/{\cal G}}}}      
\newcommand{\agbfact}[1][]{\text{$\agb/\!\aeqrel$}}    
\newcommand{\Gb}{{\overline{{\cal G}}}}                

\newcommand{\qa}{{\quer{A}}}                           

\newcommand{\holgr}{{\mathbf H}}                       
\newcommand{\bz}{{\mathbf B}}                          

\newcommand{\LG}{{\mathbf{G}}}                         
\newcommand{\aeqrelzush}[1][]{\sim}                    
\newcommand{\alg}{\mathfrak{A}}                          
\newcommand{\blg}{\mathfrak{B}}                          
\newcommand{\ideal}{\mathfrak{I}}                        
\newcommand{\malg}{{\cal M}(\alg)}                     

\newcommand{\nklza}[1][]{\ifthenelse{\equal{#1}{}}     
                                    {\rnkl{Z(\holgr_\qa)}{\LG}}        
                                   {\rnkl{Z(\holgr_{#1})}{\LG}}}       
\newcommand{\nkla}[1][]{\ifthenelse{\equal{#1}{}}      
                                    {\rnkl{\bz(\qa)}{\Gb}}        
                                    {\rnkl{\bz(#1)}{\Gb}}}       

\newcommand{\charakt}{\chi}                            

\newcommand{\YM}{{\text{YM}}}                          

\newcommand{\ymwirk}[1][]{\ifthenelse{\equal{#1}{}}{S_{\YM}}{S_{\YM,#1}}}

\newcommand{\bmat}{\begin{pmatrix}}
\newcommand{\emat}{\end{pmatrix}}

\newcommand{\ListNullAbstaende}{\setlength{\topsep}{0pt}%
                                \setlength{\parskip}{0pt}%
                                \setlength{\partopsep}{0pt}%
                                \setlength{\itemsep}{0pt}%
                                \setlength{\parsep}{0pt}}
\newcommand{\ListNurAnstrichAbstand}{\setlength{\topsep}{0pt}%
                                     \setlength{\parskip}{0pt}%
                                     \setlength{\partopsep}{0pt}%
                                     \setlength{\parsep}{0pt}}
\newenvironment{StandardListe}[2]%
               {\begin{list}%
                      {#1}%
                      {\settowidth{\leftmargin}{M#1}%
                       \settowidth{\labelwidth}{#1}%
                       \settowidth{\labelsep}{M}%
                       #2%
                      }%
                }%
               {\end{list}}%
\newenvironment{EinfachListe}[1]%
               {\begin{StandardListe}{#1}{\ListNullAbstaende}}%
               {\end{StandardListe}}%
               {\begin{StandardListe}{#1}{\ListNurAnstrichAbstand}}%
               {\end{StandardListe}}%
\newcommand{\labelsatz}[1]{#1}
\newcounter{listennr}                      %
\newlength{\hilfslaenge}
\newlength{\stdlabellaenge}
\newlength{\maximum}
\newcommand{\stdlabel}{}
\newcommand{\Maximum}{}
\newcommand{\iitem}[1][]{\ifthenelse{\equal{#1}{}}%
                           {\item \setlength{\hilfslaenge}{\stdlabellaenge}}%
                           {\item[\labelsatz{#1}\hfill]%
                            \settowidth{\hilfslaenge}{\labelsatz{#1}}}%
                         \ifthenelse{\lengthtest{\maximum < \hilfslaenge}}%
                           {\setlength{\maximum}{\hilfslaenge}%
                            \ifthenelse{\equal{#1}{}}%
                               {\renewcommand{\Maximum}{\stdlabel}}%
                               {\renewcommand{\Maximum}{#1}}}%
                           {}%
                      }      
\makeatletter
\newenvironment{AutoLabelLaengenListe}[2][]%
               {\begin{list}%
                      {\labelsatz{#1}\hfill}%
                      {\stepcounter{listennr}%
                       \settowidth{\leftmargin}{M\labelsatz{\ref{listnr\arabic{listennr}}}}%
                       \settowidth{\labelwidth}{\labelsatz{\ref{listnr\arabic{listennr}}}}%
                       \settowidth{\labelsep}{M}%
                       \settowidth{\stdlabellaenge}{\labelsatz{#1}}%
                       \renewcommand{\stdlabel}{#1}%
                       #2%
                       \renewcommand{\Maximum}{}%
                      }%
                }%
               {\renewcommand{\@currentlabel}{\Maximum}%
                \label{listnr\arabic{listennr}}%
                \end{list}%
                }%
\makeatother
\newenvironment{StandardEinrueckung}[2]%
               {\begin{list}%
                      {#1}%
                      {\settowidth{\leftmargin}{M#1}%
                       \settowidth{\labelwidth}{#1}%
                       \settowidth{\labelsep}{M}%
                       #2%
                      }%
                \item}%
               {\end{list}}%
\newenvironment{Einrueckungpur}[1]%
               {\begin{StandardEinrueckung}{#1}{\ListNullAbstaende}}%
               {\end{StandardEinrueckung}}%
\newenvironment{Einrueckung}[1]%
               {\begin{StandardEinrueckung}{#1}{\setlength{\parsep}{0pt}}}%
               {\end{StandardEinrueckung}}%

\makeatletter
\newcommand{\EineNumZeileGleichung}[2][0.5ex]
           {
            
            \vspace{#1} 
            \noindent
            \stepcounter{equation}
            \renewcommand{\@currentlabel}{\arabic{equation}}%
            \phantom{(\arabic{equation})}\hspace*{\fill}
            $\displaystyle{#2}$
            \hspace*{\fill}
            (\arabic{equation})

            \vspace{#1} 
            
           }
\makeatother
\makeatletter
\newcommand{\EineErwNumZeileGleichung}[2][0.5ex]
           {
            
            \vspace{#1} 
            \noindent
            \stepcounter{equation}
            \renewcommand{\@currentlabel}{\arabic{equation}}%
            \phantom{(\arabic{equation})}\hspace*{\fill}
            #2 %
            \hspace*{\fill}
            (\arabic{equation})

            \vspace{#1} 
            
           }
\makeatother
\newcommand{\breitrel}[1]{\hspace*{\tabcolsep} #1 \hspace*{\tabcolsep}}
\newlength{\abstaug}              %
\newenvironment{AllgUnnumGleichung}[2][1.0ex]
               {
  
                \setlength{\abstaug}{#1}
                \vspace{\abstaug}
                \hspace*{\fill}
                $\begin{array}[t]{#2}
                }%
               {\end{array}$
                \hspace*{\fill}
  
                \vspace{\abstaug}

                }%
\newenvironment{AllgNumGleichung}[2][0.0ex]
               {
  
                \setlength{\abstaug}{#1}
                \vspace{\abstaug}
                $\begin{tabular*}{\textwidth}[t]{#2}
                }%
               {\end{tabular*}$

                \vspace{\abstaug}

               }%
\newenvironment{StandardUnnumGleichungKlein}[1][0ex]
               {\renewcommand{\s}{\\[#1] }%
                \begin{AllgUnnumGleichung}{rcl}}%
               {\end{AllgUnnumGleichung}}%
\newenvironment{StandardUnnumGleichung}[1][0ex]%
               {\renewcommand{\s}{\\[#1] }%
                \begin{AllgUnnumGleichung}{>{\displaystyle}rc>{\displaystyle}l}}%
               {\end{AllgUnnumGleichung}}%
\newenvironment{XrelYZNumGleichung}[1][0ex]
               {\renewcommand{\s}{\\[#1] }%
                \begin{AllgNumGleichung}{rcll}}%
               {\end{AllgNumGleichung}}%

\newcommand{\erllang}[2][0.5\textwidth]%
              {\hfill\hspace*{1.5em}%
               \begin{minipage}[t]{#1}{\small%
                          \begin{list}{(}{\ListNullAbstaende%
                                          \settowidth{\leftmargin}{(}%
                                          \settowidth{\labelwidth}{(}%
                                          \settowidth{\labelsep}{}%
                                         }%
                          \item#2)%
                          \end{list}}%
               \end{minipage}\\[-0.9ex]
              }%
\newcommand{\DefBemUmgeb}[1]%
           {\newenvironment{#1}[1][]%
                           {\begin{Einrueckung}{{\bf #1}}%
                            \ifx##1\empty\else{{\bf ##1}
                            
                                                        }\fi%
                            }%
                           {\end{Einrueckung}}}
\newcommand{\DefSBemUmgeb}[2]
           {\newenvironment{#1}[1][]%
                           {\begin{Einrueckung}{{\bf #2}}%
                            \ifx##1\empty\else{{\bf ##1}
                            
                                                        }\fi%
                            }%
                           {\end{Einrueckung}}}
\makeatletter
\newcommand{\DefBspUmgeb}[3]
           {\newcounter{#2}[#3]%
            \newenvironment{#1}[1][]%
                           {\stepcounter{#2}%
                            \renewcommand{\ZaehlerMarke}{\arabic{#2}}%
                            \renewcommand{\Einzugsname}{{\bf #1 \ZaehlerMarke}}%
                            \begin{Einrueckung}{\Einzugsname}
                            \ifx##1\empty\else{{\bf ##1}\\}\fi%
                            \renewcommand{\@currentlabel}{\ZaehlerMarke}%
                            }%
                           {\end{Einrueckung}}}
\makeatother
\newcommand{\ZaehlerbisEbene}{section}
\newcommand{\Ebenea}{section}
\newcommand{\Ebeneb}{subsection}

\newcommand{\Abschnittnummer}{%
            \ifx\ZaehlerbisEbene\Ebenea{\arabic{section}}%
             \else{%
              \ifx\ZaehlerbisEbene\Ebeneb{\arabic{section}.\arabic{subsection}}%
               \else{\arabic{section}.\arabic{subsection}.\arabic{subsubsection}}%
              \fi}%
            \fi}     
\newcommand{\Einzugsname}{}
\newcommand{\ZaehlerMarke}{}
\makeatletter
\newcommand{\DefThmUmgeb}[3]%
           {\newcounter{#1}[#3]%
            \newenvironment{#1}[1][]%
                           {\stepcounter{#2}%
                            \setcounter{#1}{\value{#2}}%
                            \renewcommand{\ZaehlerMarke}{\Abschnittnummerpunkt\arabic{#1}}%
                            \renewcommand{\Einzugsname}{{\bf #1 \ZaehlerMarke}}%
                            \begin{Einrueckung}{\Einzugsname}
                            \ifx##1\empty\else{{\bf ##1}
                            
                                                        }\fi%
                            \renewcommand{\@currentlabel}{\ZaehlerMarke}%
                            }%
                           {\end{Einrueckung}}}
\makeatother
\makeatletter
\newcommand{\DefSThmUmgeb}[4]%
           {\newcounter{#1}[#3]%
            \newenvironment{#1}[1][]%
                           {\stepcounter{#2}%
                            \setcounter{#1}{\value{#2}}%
                            \renewcommand{\ZaehlerMarke}{\Abschnittnummerpunkt\arabic{#1}}%
                            \renewcommand{\Einzugsname}{{\bf #4 \ZaehlerMarke}}
                            \begin{Einrueckung}{\Einzugsname}
                            \ifx##1\empty\else{{\bf ##1}

                                                        }\fi%
                            \renewcommand{\@currentlabel}{\ZaehlerMarke}%
                            }%
                           {\end{Einrueckung}}}
\makeatother
\makeatletter
\newcommand{\DefUnterNumThmUmgeb}[5]%
           {\newcounter{#1}[#3]%
            \newcounter{#4}%
            \newenvironment{#1}[1][]%
                           {\ifx##1\empty\else{\stepcounter{#2}\setcounter{#4}{0}}\fi%
                            \stepcounter{#4}%
                            \setcounter{#1}{\value{#2}}%
                            \renewcommand{\ZaehlerMarke}{\Abschnittnummerpunkt\arabic{#1}\alph{#4}}%
                            \renewcommand{\Einzugsname}{{\bf #5 \ZaehlerMarke}}
                            \begin{Einrueckung}{\Einzugsname}
                            \renewcommand{\@currentlabel}{\ZaehlerMarke}%
                            }%
                           {\end{Einrueckung}}}
\makeatother
\newenvironment{Beweis}[1][]%
               {\begin{Einrueckung}{{\bf Beweis}}%
                \ifx#1\empty\else{{\bf #1}

                                            }\fi%
                }%
               {\end{Einrueckung}%
                }%
\newenvironment{Proof}[1][]%
               {\begin{Einrueckung}{{\bf Proof}}%
                \ifx#1\empty\else{{\bf #1}

                                            }\fi%
                }%
               {\end{Einrueckung}%
                }%
               {\begin{Einrueckung}{{\bf \glqq Beweis\grqq}}%
                \ifx#1\empty\else{{\bf #1}
                
                                            }\fi%
                }%
               {\end{Einrueckung}%
                }%
               {\begin{Einrueckung}{{\bf Begr"undung}}%
                \ifx#1\empty\else{{\bf #1}
                
                                            }\fi%
                }%
               {\end{Einrueckung}%
                }%
\newenvironment{Hinrichtung}%
               {\begin{Einrueckungpur}{$\impliz$}}%
               {\end{Einrueckungpur}}%
\newenvironment{Rueckrichtung}%
               {\begin{Einrueckungpur}{$\invimpliz$}}%
               {\end{Einrueckungpur}}%
               {\begin{Einrueckungpur}{\glqq$\teilmenge$\grqq}}%
               {\end{Einrueckungpur}}%
               {\begin{Einrueckungpur}{\glqq$\obermenge$\grqq}}%
               {\end{Einrueckungpur}}%
               {\begin{Einrueckungpur}{"$\teilmenge$"}}%
               {\end{Einrueckungpur}}%
               {\begin{Einrueckungpur}{"$\obermenge$"}}%
               {\end{Einrueckungpur}}%
\newcommand{\qed}{\nopagebreak\hspace*{2em}\hspace*{\fill}{\bf qed}}
\newcommand{\ARabic}{\arabic}
\newcommand{\Nummerntypa}{\arabic}   
\newcommand{\Nummerntypb}{\alph}
\newcommand{\Nummerntypc}{\roman}
\newcommand{\Nummerntypd}{\Alph}

\newcommand{\Nra}{\Nummerntypa{Nummera}}            %
\newcommand{\Nrb}{\Nummerntypb{Nummerb}}            %
\newcommand{\Nrc}{\Nummerntypc{Nummerc}}                
\newcommand{\Nrd}{\Nummerntypd{Nummerd}}                
\newcommand{\ZeichenzuNrTyp}[1]%
           {\ifx#1\ARabic {.}\else{)}%
                  \fi}                              %
\newcommand{\NrZeicha}{\ZeichenzuNrTyp{\Nummerntypa}}
\newcommand{\NrZeichb}{\ZeichenzuNrTyp{\Nummerntypb}}
\newcommand{\NrZeichc}{\ZeichenzuNrTyp{\Nummerntypc}}
\newcommand{\NrZeichd}{\ZeichenzuNrTyp{\Nummerntypd}}
\newcommand{\ListMarkea}%
           {\Nra\NrZeicha}
\newcommand{\ListMarkeb}%
           {\Nra\NrZeicha\Nrb\NrZeichb}
\newcommand{\ListMarkec}%
           {\Nra\NrZeicha\Nrb\NrZeichb\Nrc\NrZeichc}
\newcommand{\ListMarked}%
           {\Nra\NrZeicha\Nrb\NrZeichb\Nrc\NrZeichc\Nrd\NrZeichd}
\newcommand{\Anfangszeichen}{}
\newcommand{\Anfangspunkt}{}
\newcounter{Schachtelebene}
\newcounter{Hilfszaehler}
\newcommand{\Hilfsbefehl}{}
\newcommand{\Schachtelebene}{\alph{Schachtelebene}}
\makeatletter
\newenvironment{AllgNumerierteListe}[2][]
               {\addtocounter{Schachtelebene}{1}%
		\setcounter{Hilfszaehler}{#2}%
                \renewcommand{\Anfangszeichen}%
                             {\renewcommand{\Hilfsbefehl}{\csname Nummerntyp\Schachtelebene \endcsname}%
                              \Hilfsbefehl{Hilfszaehler}}%
                \renewcommand{\Anfangspunkt}%
                             {\csname NrZeich\Schachtelebene \endcsname}%
                \begin{list}%
                      {\stepcounter{Nummer\Schachtelebene}%
                       \csname Nr\Schachtelebene \endcsname
                       \csname NrZeich\Schachtelebene \endcsname
                       }%
                      {\settowidth{\leftmargin}{M\Anfangszeichen\Anfangspunkt}%
                       \settowidth{\labelwidth}{\Anfangszeichen\Anfangspunkt}%
                       \settowidth{\labelsep}{M}%
                       \setlength{\topsep}{0pt}%
                       \setlength{\parskip}{0pt}%
                       \setlength{\partopsep}{0pt}%
                       \setlength{\itemsep}{0pt}%
                       \setlength{\parsep}{0pt}%
                      }%
                \renewcommand{\@currentlabel}{\csname ListMarke\Schachtelebene \endcsname}%
                }%
               {\ifthenelse{\equal{}{}}{\setcounter{Nummer\Schachtelebene}{0}}{}
                \addtocounter{Schachtelebene}{-1}%
                \end{list}}
\makeatother
\newenvironment{NumerierteListe}[1]
               {\begin{AllgNumerierteListe}{#1}}
               {\end{AllgNumerierteListe}}
\newenvironment{WeiterNumerierteListe}[1]
               {\begin{AllgNumerierteListe}[Weiter]{#1}}
               {\end{AllgNumerierteListe}}

\newcommand{\UnnumAnfangszeichen}{}
\newcounter{UnnumSchachtelebene}
\newcommand{\UnnumSchachtelebene}{\alph{UnnumSchachtelebene}}
\makeatletter
\newenvironment{UnnumerierteListe}%
               {\addtocounter{UnnumSchachtelebene}{1}%
                \renewcommand{\UnnumAnfangszeichen}%
                             {\csname UnnumZeich\UnnumSchachtelebene \endcsname}%
                \begin{list}%
                      {\UnnumAnfangszeichen}%
                      {\settowidth{\leftmargin}{M\UnnumAnfangszeichen}%
                       \settowidth{\labelwidth}{\UnnumAnfangszeichen}%
                       \settowidth{\labelsep}{M}%
                       \setlength{\topsep}{0pt}%
                       \setlength{\parskip}{0pt}%
                       \setlength{\partopsep}{0pt}%
                       \setlength{\itemsep}{0pt}%
                       \setlength{\parsep}{0pt}%
                      }%
                }%
               {\addtocounter{UnnumSchachtelebene}{-1}%
                \end{list}}
\makeatother
\newlength{\fktdefhilfslaenge}
\newcommand{\ohnefktdef}[4]
           {\hspace*{\fill}
            $\begin{array}[t]{ccc}%
            #1 & \nach & #2 \\
            #3 & \auf  & #4
            \end{array}$
            \hspace*{\fill}}
\newcommand{\fktdef}[5]
           {\hspace*{\fill}
            $\begin{array}[t]{cccc}%
            #1: & #2 & \nach & #3 \\    
                & #4 & \auf  & #5
            \end{array}$
            \settowidth{\fktdefhilfslaenge}{$#1$:}
            \hspace*{0.6 \fktdefhilfslaenge}  
            \hspace*{\fill}}
\newcommand{\fktdefpur}[5]
           {$\begin{array}[t]{cccc}%
            #1: & #2 & \nach & #3 \\    
                & #4 & \auf  & #5
            \end{array}$}
\newcommand{\fktdefabgesetztpur}[5]
           {
            
            $\begin{array}[t]{cccc}%
            #1: & #2 & \nach & #3 \\    
                & #4 & \auf  & #5
            \end{array}$
            \settowidth{\fktdefhilfslaenge}{$#1$:}
            \hspace*{0.6 \fktdefhilfslaenge}
            
           }
\newcommand{\fktdefabgesetzt}[5]
           {
           
            \hspace*{\fill}
            $\begin{array}[t]{cccc}%
            #1: & #2 & \nach & #3 \\    
                & #4 & \auf  & #5
            \end{array}$
            \settowidth{\fktdefhilfslaenge}{$#1$:}
            \hspace*{0.6 \fktdefhilfslaenge}  
            \hspace*{\fill}
            
            }
\newcommand{\ohnefktdefabgesetzt}[4]
           {      

            \hspace*{\fill}
            $\begin{array}[t]{ccc}%
            #1 & \nach & #2 \\
            #3 & \auf  & #4
            \end{array}$
            \hspace*{\fill}

            }
\newcommand{\doppelohnefktdefabgesetzt}[6]
           {

            \hspace*{\fill}
            $\begin{array}[t]{ccccc}%
            #1 & \nach & #2 & \nach & #3\\
            #4 & \auf  & #5 & \auf  & #6
            \end{array}$
            \hspace*{\fill}

            }
\newcommand{\anhang}%
           {\appendix
            \sectioninh{Anhang}
            \renewcommand{\Abschnittnummer}{%
                  \ifx\ZaehlerbisEbene\Ebenea{\Alph{section}}%
                  \else{%
                        \ifx\ZaehlerbisEbene\Ebeneb{\Alph{section}.\arabic{subsection}}%
                        \else{\Alph{section}.\arabic{subsection}.\arabic{subsubsection}}%
                        \fi}%
                  \fi}%
                 
            }            
\newcommand{\anhangengl}%
           {\appendix
            \sectioninh{Appendix}
            \renewcommand{\Abschnittnummer}{%
                  \ifx\ZaehlerbisEbene\Ebenea{\Alph{section}}%
                  \else{%
                        \ifx\ZaehlerbisEbene\Ebeneb{\Alph{section}.\arabic{subsection}}%
                        \else{\Alph{section}.\arabic{subsection}.\arabic{subsubsection}}%
                        \fi}%
                  \fi}%
                 
            }

\newcounter{wdhlstufe}
\newcommand{\sectioninh}[1]%
           {\section*{#1}%
            \addcontentsline{toc}{section}{#1}}
\newcommand{\bezeichnung}[3]%
           {\begin{Einrueckungpur}{\hbox to 6em{#1}\hbox to 2.4em{\hfill#2}}
            #3
            \end{Einrueckungpur}}

\newcommand{\doppelteinfach}{e}

\newcommand{\ifdoppelt}[1]{\ifthenelse{\equal{\doppelteinfach}{d}}{#1}{}}
\newcommand{\ifeinfach}[1]{\ifthenelse{\equal{\doppelteinfach}{e}}{#1}{}}

\newlength{\querfhilfsl}              %

\newlength{\hll}

\DefThmUmgeb{Theorem}{Theorem}{\ZaehlerbisEbene}
\DefThmUmgeb{Definition}{Definition}{\ZaehlerbisEbene}
\DefThmUmgeb{Satz}{Theorem}{\ZaehlerbisEbene}
\DefThmUmgeb{Proposition}{Theorem}{\ZaehlerbisEbene}
\DefThmUmgeb{Lemma}{Theorem}{\ZaehlerbisEbene}
\DefThmUmgeb{Folgerung}{Theorem}{\ZaehlerbisEbene}
\DefThmUmgeb{Corollary}{Theorem}{\ZaehlerbisEbene}
\DefThmUmgeb{Vorschrift}{Definition}{\ZaehlerbisEbene}
\DefThmUmgeb{Construction}{Definition}{\ZaehlerbisEbene}
\DefSThmUmgeb{FormSatz}{Theorem}{\ZaehlerbisEbene}{\glqq Satz\grqq} 
\DefThmUmgeb{Vermutung}{Theorem}{\ZaehlerbisEbene}
\DefThmUmgeb{Conjecture}{Theorem}{\ZaehlerbisEbene}
\DefThmUmgeb{Konvention}{Definition}{\ZaehlerbisEbene}
\DefThmUmgeb{Feststellung}{Theorem}{\ZaehlerbisEbene}
\DefUnterNumThmUmgeb{DefinitionZusatzNum}{Definition}{\ZaehlerbisEbene}{DefZN}{Definition}
\DefBspUmgeb{Beispiel}{Beispiel}{subsubsection}
\DefBspUmgeb{Example}{Example}{subsubsection}
\DefBspUmgeb{Frage}{Frage}{section}
\DefBspUmgeb{Question}{Question}{section}
\DefBspUmgeb{Aufgabe}{Aufgabe}{section}
\DefBspUmgeb{Ziel}{Ziel}{section}
\DefBemUmgeb{Bemerkung}
\DefBemUmgeb{Remark}
\DefSBemUmgeb{OffeneFrage}{Offene Frage}
\newcommand{\bdf}{\begin{Definition}}
\newcommand{\edf}{\end{Definition}}
\newcommand{\bvorsch}{\begin{Vorschrift}}
\newcommand{\evorsch}{\end{Vorschrift}}
\newcommand{\bconst}{\begin{Construction}}
\newcommand{\econst}{\end{Construction}}
\newcommand{\bthm}{\begin{Theorem}}
\newcommand{\ethm}{\end{Theorem}}
\newcommand{\bsatz}{\begin{Satz}}
\newcommand{\esatz}{\end{Satz}}
\newcommand{\bprop}{\begin{Proposition}}
\newcommand{\eprop}{\end{Proposition}}
\newcommand{\blem}{\begin{Lemma}}
\newcommand{\elem}{\end{Lemma}}
\newcommand{\bfolg}{\begin{Folgerung}}
\newcommand{\efolg}{\end{Folgerung}}
\newcommand{\bcorr}{\begin{Corollary}}
\newcommand{\ecorr}{\end{Corollary}}
\newcommand{\bfest}{\begin{Feststellung}}
\newcommand{\efest}{\end{Feststellung}}
\newcommand{\bbew}{\begin{Beweis}}
\newcommand{\ebew}{\end{Beweis}}
\newcommand{\bpf}{\begin{Proof}}
\newcommand{\epf}{\end{Proof}}
\newcommand{\bwnum}{\begin{WeiterNumerierteListe}}
\newcommand{\ewnum}{\end{WeiterNumerierteListe}}
\newcommand{\bdfzn}{\begin{DefinitionZusatzNum}}
\newcommand{\edfzn}{\end{DefinitionZusatzNum}}
\newcommand{\bbem}{\begin{Bemerkung}}
\newcommand{\ebem}{\end{Bemerkung}}
\newcommand{\brem}{\begin{Remark}}
\newcommand{\erem}{\end{Remark}}
\newcommand{\bnum}{\begin{NumerierteListe}}
\newcommand{\enum}{\end{NumerierteListe}}
\newcommand{\bunum}{\begin{UnnumerierteListe}}
\newcommand{\eunum}{\end{UnnumerierteListe}}
\newcommand{\bbsp}{\begin{Beispiel}}
\newcommand{\ebsp}{\end{Beispiel}}
\newcommand{\bex}{\begin{Example}}
\newcommand{\eex}{\end{Example}}
\newcommand{\bfrag}{\begin{Frage}}
\newcommand{\efrag}{\end{Frage}}
\newcommand{\bquest}{\begin{Question}}
\newcommand{\equest}{\end{Question}}
\newcommand{\baufg}{\begin{Aufgabe}}
\newcommand{\eaufg}{\end{Aufgabe}}
\newcommand{\bof}{\begin{OffeneFrage}}
\newcommand{\eof}{\end{OffeneFrage}}
\newcommand{\bverm}{\begin{Vermutung}}
\newcommand{\everm}{\end{Vermutung}}
\newcommand{\bconj}{\begin{Conjecture}}
\newcommand{\econj}{\end{Conjecture}}
\newcommand{\bkonv}{\begin{Konvention}}
\newcommand{\ekonv}{\end{Konvention}}
\newcommand{\bglklein}{\begin{StandardUnnumGleichungKlein}}
\newcommand{\eglklein}{\end{StandardUnnumGleichungKlein}}
\newcommand{\bgl}{\begin{StandardUnnumGleichung}}
\newcommand{\egl}{\end{StandardUnnumGleichung}}
\newcommand{\bglrtext}{\begin{XrelYZNumGleichung}}
\newcommand{\eglrtext}{\end{XrelYZNumGleichung}}

\newcommand{\berlgl}{\begin{StandardUnnumGleichung}}
\newcommand{\eerlgl}{\end{StandardUnnumGleichung}}
\newcommand{\beinrueck}{\begin{Einrueckungpur}} 
\newcommand{\eeinrueck}{\end{Einrueckungpur}}
\newcommand{\beinflist}{\begin{EinfachListe}} 
\newcommand{\eeinflist}{\end{EinfachListe}}
\newcommand{\beq}{\begin{equation}}
\newcommand{\eeq}{\end{equation}}
\newcommand{\bhin}{\begin{Hinrichtung}}
\newcommand{\ehin}{\end{Hinrichtung}}
\newcommand{\brueck}{\begin{Rueckrichtung}}
\newcommand{\erueck}{\end{Rueckrichtung}}
\newcommand{\bvl}{\begin{AutoLabelLaengenListe}{\ListNullAbstaende}}
\newcommand{\evl}{\end{AutoLabelLaengenListe}}
\newcommand{\df}[1]{{\bf #1}}

\newlength{\adressabstand}
\usepackage{diagrams}
\newcommand{\borelalg}{{\cal B}}

\newcommand{\twist}{f}
\newcommand{\alga}{\alg_0}
\newcommand{\algb}{\alg_1}
\newcommand{\algab}{\alg}
\newcommand{\algablang}{\sumalglang}
\newcommand{\malga}{\malg_0}
\newcommand{\malgb}{\malg_1}
\newcommand{\malgab}{\malg}

\newcommand{\elalga}{\elalg_0}
\newcommand{\elalgb}{\elalg_1}
\newcommand{\elalgab}{\elalg}

\newcommand{\iotaalgb}{\iota_1}
\newcommand{\iotaalgab}{\iota}

\newcommand{\charaktalgb}{\varphi}

\newcommand{\offenCalga}{\offenCalg_0}
\newcommand{\offenCalgb}{\offenCalg_1}

\newcommand{\offenCalgboth}{\offenCalga, \offenCalgb} 
\renewcommand{\offenCalgboth}{U} 
\renewcommand{\offenCalga}{U} 
\renewcommand{\offenCalgb}{U} 
 
\renewcommand{\alg}{{\mathfrak A}}
\renewcommand{\blg}{{\mathfrak B}}
\newcommand{\clg}{{\mathfrak C}}

\newcommand{\mlg}{\mathfrak{M}}
\newcommand{\allgsumalg}{\mlg + \alg}

\newcommand{\elalg}{a} 
\newcommand{\elblg}{b} 
\newcommand{\elclg}{c} 
 
\newcommand{\elmlg}{m}
\newcommand{\elideal}{i}
\renewcommand{\elideal}{n}
\renewcommand{\mlg}{\ideal}
\renewcommand{\elmlg}{\elideal}

\newcommand{\offenCalg}{A}

\renewcommand{\malg}{\spec \alg}

\newcommand{\sumalglang}{\alga \dirvsum \algb}
\newcommand{\msumalglang}{\spec (\sumalglang)}

\newcommand{\cnulleinschr}[2]{C_{0,#1}(#2)}
\newcommand{\elsubcover}{U}

\newcommand{\set}{\mathbf S}
\newcommand{\elset}{s} 
\newcommand{\Y}{\mathbf{Y}} 
\renewcommand{\Y}{Y} 
\newcommand{\ely}{\mathbf{y}} 
\renewcommand{\ely}{y}

\newcommand{\charabb}{\tau}
\newcommand{\charact}{\charakt}
\newcommand{\disjunion}{\sqcup}
\newcommand{\setabb}{\sigma}

\newcommand{\cover}{{\cal O}}
\newcommand{\subcover}{{\cal U}}

\newcommand{\gelf}[1][]{\ifthenelse{\equal{#1}{}}{\widetilde}{G_{#1}}}

\renewcommand{\set}{X}
\renewcommand{\elset}{x}

\newcommand{\nonicht}[1]{\textbf{Das ist hier noch unbesetzt.}}

\newcommand{\AP}{{\mathrm{AP}}}

\newarrow{auf}|---{->}
\newarrow{nach}----{->}
\newarrow{ident}33333
\newarrow{aufsurj}|---{->>}
\newarrow{nachsurj}----{->>}
\newarrow{einbett}C---{->}
\newarrow{strichel}{}{dash}{}{dash}{->}
\newarrow{impliz}===={=>}
\newarrow{fillauf}....{->}
\newlength{\CDhoehe}                  
\newlength{\CDgap}                    
\makeatletter
\newcounter{tab}
\newcommand{\fusszeile}{}
                        {\par\noindent\fusszeile
                         \end{center}}
\makeatother
\DefThmUmgeb{Convention}{Definition}{\ZaehlerbisEbene}
\DefThmUmgeb{Notation}{Definition}{\ZaehlerbisEbene}

\renewcommand{\boundfkt}{C_b}

\newcommand{\Bigbetrag}[1]%
           {\Bigl|{#1}\Bigr|}            
\newcommand{\Biggbetrag}[1]%
           {\Biggl|{#1}\Biggr|}

\newcommand{\Bohr}{{\mathrm{Bohr}}}

\newcommand{\voffen}{{V}}
\newcommand{\vkomp}{{K}}
\newcommand{\compl}[1]{\complement#1}

\newcommand{\bconv}{\begin{Convention}}
\newcommand{\econv}{\end{Convention}}
\newcommand{\bnotat}{\begin{Notation}}
\newcommand{\enotat}{\end{Notation}}

\newcommand{\bpm}{\begin{pmatrix}}
\newcommand{\epm}{\end{pmatrix}}

\DeclareMathOperator{\spec}{spec}

\DeclareMathOperator{\dirvsum}{\boxplus}
\renewcommand{\dirvsum}{\boxplus}

\newcommand{\bnotatlist}{\bvl{}}
\newcommand{\enotatlist}{\evl{}}
\newcommand{\citem}[1][]{\iitem[$\bullet$~~ $#1$]$\ldots$\quad\!\!}
\begin{document}
\title{Spectra of Abelian $C^\ast$-Subalgebra Sums}
\author{Christian Fleischhack\thanks{e-mail: 
            {\tt fleischh@math.upb.de}} \\   
        \\
        {\normalsize\em Institut f\"ur Mathematik}\\[\adressabstand]
        {\normalsize\em Universit\"at Paderborn}\\[\adressabstand]
        {\normalsize\em Warburger Stra\ss e 100}\\[\adressabstand]
        {\normalsize\em 33098 Paderborn}\\[\adressabstand]
        {\normalsize\em Germany}
        \\[-25\adressabstand]}     
\date{September 18, 2014}
\maketitle
\newcommand{\iotarestr}{\iota_\stdrestr}
\newcommand{\stdrestr}{\setabb}

\newcommand{\algxxxnull}{\alg_0}
\newcommand{\algxxxeins}{\alg_1}
\newcommand{\malgxxxnull}{\malg_0}
\newcommand{\malgxxxeins}{\malg_1}
\newcommand{\elalgxxxnull}{\elalg_0}
\newcommand{\elalgxxxeins}{\elalg_1}

\begin{abstract}
Let 
$\boundfkt(\set)$ be the $C^\ast$-algebra of bounded continuous functions on 
some non-compact, but locally compact Hausdorff space $X$.
Moreover, let $\alga$ be some ideal 
and $\algb$ be some unital $C^\ast$-sub\-al\-ge\-bra
of $\boundfkt(\set)$.
For $\alga$ and $\algb$ having trivial intersection, we show
that the spectrum of their vector space sum
equals the disjoint union of their individual spectra,
whereas their topologies are nontrivially interwoven. Indeed, they
form a so-called twisted-sum topology which we will investigate before.
Within the whole framework, e.g., the one-point compactification of $X$
and the spectrum of the algebra of asymptotically almost periodic functions
can be described.

\end{abstract}

\section{Introduction}
It is well known that the spectrum of a direct sum of two abelian $C^\ast$-algebras 
equals the topological direct sum of the respective individual spectra.
Sometimes, however, one is given only a {\em vector space}\/ direct
sum of two $C^\ast$-algebras.
This applies, most prominently, to $\alg + \C \EINS$ to be considered when one 
adjoins a unit to 
the non-unital $C^\ast$-algebra $\alg$. 
Another example is the algebra $C_0(X) + C_\AP(X)$ of 
asymptotically almost periodic functions \cite{GrigTonev} 
on a non-compact locally compact abelian group $X$, being our
main motivation \cite{paper39}.
It is now natural to ask whether there 
are still general arguments on how to determine the spectrum in these cases.
Or, to put it into a more abstract form:
how does the spectrum of 
a sum
of any two abelian $C^\ast$-algebras look like?

\enlargethispage{0.5\baselineskip}

Of course, in this generality, the question makes no sense, as we do not know
how to multiply elements of different addends.
Even if both algebras are contained in a third algebra, their sum need
not form an algebra again. Therefore, let us take a closer look at the 
situations above. In both cases, 
we are given two $C^\ast$-algebras 
$\alga$ and $\algb$ that fulfill 
$\alga \algb \teilmenge \alga$ and that trivially intersect. Here, 
their direct vector space sum $\sumalglang$ is, at least, a $\ast$-algebra.
However, in order to get a $C^\ast$-algebra structure, we need 
a norm. In both cases above, this does not cause a problem as 
$\alg_0$ and $\alg_1$ are
contained in some $C^\ast$-algebra $\clg$. We shall assume this in the following,
as then, by general arguments 
(Corollary \ref{corr:quasiideal+subalg=Cstern}),
$\alga + \algb$ is a $C^\ast$-algebra containing $\alga$ as an ideal.
Note that, since $C^\ast$-norms are unique on $\ast$-algebras, the
whole construction is independent of the choice of $\clg$.

As we are going to calculate spectra of \emph{abelian} $C^\ast$-algebras,
we may assume
that $\clg$ 
equals the $C^\ast$-algebra $C_0(X)$ of continuous vanishing-at-infinity functions 
on some locally compact $X$.
Assuming for the moment that $\alga$ is an ideal also in $\clg$, it is 
necessarily \cite{ConwayFktA} 
of the form
\bgl
\alga \breitrel\iso C_{0,\Y}(X) 
 & := & \bigl\{f \in C_0(X) \mid \text{$f \ident 0$ on $\compl \Y$}\bigr\}
 \breitrel\iso C_0(\Y)
\egl\noindent
for some open $\Y \teilmenge X$. 
This, however, is not the perfect
framework for the case we are interested in most,
namely the asymptotically almost periodic functions.
Here, one is tempted to 
choose $\Y = X = \R$, but then $C_\AP(\R)$ is not a
subalgebra of $C_0(\R)$. 
As we are aiming at unital $\algb$ anyway and since we 
are free to choose any $\clg$ containing $\alga$ and $\algb$,
we will therefore prefer $\clg$ to be the set $\boundfkt(X)$ of all 
bounded continuous functions on $X$.
To wrap up, we will let $\alga$ be an ideal and $\algb$ be 
a unital subalgebra of $\boundfkt(X)$ with $\alga \cap \algb = \NULL$.

What can one say now about 
the spectrum of the sum $\alga + \algb$ of these two algebras?
Considering the unitization, i.e., $\alga = C_0(X)$ and $\algb = \C \, \EINS$
as subalgebras of $\boundfkt(X)$, we see that, as a set, the spectrum of the sum is
the disjoint union of the single spectra, namely $X$ and $\{\infty\}$.
However, there are certain matching conditions influencing the topology.
In fact, the topology is \emph{not} generated by the open sets in $X$ and
in $\{\infty\}$; it is given by the open sets in $X$ and by complements of closed compacta
in $X$ together with $\infty$. 
In other words, the spectrum of the sum is the disjoint union of the two spectra,
but their topologies get nontrivially interwoven.
This will indeed remain true for the general case.
To see this, we will construct an appropriate isomorphism
\bgl
\charabb \breitrel: \malga \disjunion \malgb \nach \msumalglang \,.
\egl\noindent
On $\malgb$, the map $\charabb$ should be given
by $[\charabb(\charaktalgb)](\elalga + \elalgb) := \charaktalgb(\elalgb)$.
Indeed, $\charabb(\charaktalgb)$ is a character on $\algb$ since
$\alga \algb \teilmenge \alga$ (check directly or see Theorem \ref{thm:spec_sum}).
On $\malga$, which we may assume to be an open 
subset $\Y$ of $X$, the situation is simpler: Here, we just set 
$[\charabb(\ely)](\elalg) := \elalg(\ely)$.
Taking an appropriate basis of the Gelfand topology on $\msumalglang$,
we will get a simple description of the topology 
on $\malga \disjunion \malgb$ as mediated by $\charabb$.
It explicitly shows how the topologies of the spectra of $\alga$ and $\algb$
are getting intertwined.
This way, in particular, we 
generalize 
the results of Grigoryan and Tonev \cite{GrigTonev} on asymptotically
almost periodic functions from $\R$ to arbitrary non-compact, but locally compact 
abelian Hausdorff groups.

\leerezeile

The paper is organized as follows: We will start in Section \ref{sect:twisted-top} 
with an abstract definition
of the so-called twisted-sum topology on the disjoint union of topological spaces.
This definition, of course, is directly motivated by the topology on 
$\msumalglang$ to be derived in Section~\ref{sect:sumspec_top}.
Before, in Section \ref{sect:c-ast-algs}
we will summarize some general facts we need from the theory of (abelian) $C^\ast$-algebras.
In Section \ref{sect:denseness}, we study how $X$ is contained in the spectrum of
$\sumalglang$. We close in Section \ref{sect:examples} with applications to 
the unitization and to asymptotically almost periodic functions.
In Appendix \ref{app:measure}, we discuss measures on the spectra.

\newcommand{\dummyraum}{M}
\newcommand{\dummyfkt}{g}
\newcommand{\dummyraumb}{N}
\newcommand{\dummyfktb}{h}
\newcommand{\rauma}{\mathbf{Y}}
\newcommand{\raumb}{\mathbf{Z}}
\newcommand{\raumc}{\mathbf{M}}
\newcommand{\raumab}{\rauma \disjunion \raumb}
\newcommand{\elrauma}{y}
\newcommand{\elraumb}{z}
\newcommand{\typzd}{23}
\newcommand{\randrest}[2][]{\ifthenelse{\equal{#1}{}}{\quer{#2}\setminus#2}{\quer{#2}\setminus(#2)}}
\newcommand{\twistabb}{f}
\newcommand{\twisttopallg}[3]{#1 \disjunion_{#3} #2}
\newcommand{\twisttopsimple}[3]{#1 \disjunion #2}
\newcommand{\stdtwisttop}{\twisttop\rauma\raumb\twistabb}
\newcommand{\twisttop}{\twisttopallg}
\newcommand{\simplifytwist}{\renewcommand{\twisttop}{\twisttopsimple}}
\newcommand{\explicittwist}{\renewcommand{\twisttop}{\twisttopallg}}

\section{Twisted Sum}
\label{sect:twisted-top}

We are first going to describe a topology on the disjoint union of two topological
spaces that is, in general, different from the standard direct sum.

\bnotat
Within this section, we let be
\bnotatlist
\citem[\rauma, \raumb] some disjoint topological spaces
\citem[\twistabb] some continuous map $\twistabb : \rauma \nach \raumb$
\citem[\compl] the complement within $\rauma$.
\enotatlist
\enotat

\bdf
\label{def:top_disjunion}
The \df{$\twistabb$-twisted topology} on the disjoint union 
$\rauma \disjunion \raumb$ is the topology generated 
by
all sets of the following types:
\bgl
\begin{array}{lcrclcl}
\text{Type 1:} && \voffen & \!\!\!\sqcup\!\!\! & \leeremenge 
        && \text{with open $\voffen \teilmenge \rauma$} \\
\text{Type 2:} && \compl\vkomp & \!\!\!\sqcup\!\!\! & \raumb 
        && \text{with compact closed $\vkomp \teilmenge \rauma$} \\
\text{Type 3:} && \twistabb^{-1}(W) & \!\!\!\sqcup\!\!\! & W
        && \text{with open $W \teilmenge \raumb$.}
\end{array}
\egl
The sets above are called \df{standard sets}.

We denote $\rauma \disjunion \raumb$, equipped with
the $\twistabb$-twisted topology, by 
\bgl
\stdtwisttop
\egl
and call it \df{$\twistabb$-twisted sum} of $\rauma$ and $\raumb$.
\edf
If we speak of \df{sum topologies} below, we will mean both the
$\twistabb$-twisted sum and the direct sum on $\rauma \disjunion \raumb$.
If $\twistabb$ is clear from the context, we may drop it.
Also note that 
$\twisttop\leeremenge\raumb\twistabb = \raumb$.

\blem
\label{lem:top_sum}
A basis for the topology on $\stdtwisttop$
is given by the following sets:
\bnum4
\item
sets of type 1;
\item
sets of type $\typzd$, i.e., intersections of a set of type 2 with a set of type 3.
\enum
\elem

\bpf
First note that the type is preserved under finite intersections of
same-type sets. Next, any intersection of a type-1 set with any standard set
is again of type 1; for this, simply observe that the $\rauma$-part 
of any standard set is open in $\rauma$.
\qed
\epf
Note that the total space $\rauma \disjunion \raumb$ is both a type-2 and
a type-3 set.

\blem
\label{lem:rel-top}
Let $\rauma \disjunion \raumb$ be given the $\twistabb$-twisted topology from 
Definition \ref{def:top_disjunion}.
Then we have:
\bunum
\item
The relative topology on $\rauma$ coincides with the original topology on $\rauma$.
\item
The relative topology on $\raumb$ coincides with the original topology on $\raumb$.
\eunum
\elem
Moreover, $\rauma$ is open and $\raumb$ is closed in $\stdtwisttop$.
\bpf
Obvious.
\qed
\epf
If confusion is unlikely, we may write $\rauma$ and $\raumb$ instead of
$\rauma \disjunion \leeremenge$ and $\leeremenge \disjunion \raumb$, respectively.

\blem
\label{lem:hausdorff-krit}
$\stdtwisttop$ is Hausdorff iff $\rauma$ is locally compact and both $\rauma$ and $\raumb$ are Hausdorff.
\elem
\bpf
\brueck
Obviously, any two distinct points in $\rauma$ can be separated by type-1 sets. Similarly,
any two distinct points in $\raumb$ can be separated by type-3 sets. Let now $\elrauma \in \rauma$ 
and $\elraumb \in \raumb$. Choose some neighbourhood $\voffen$ of $\elrauma$ 
contained in some compactum $\vkomp$.
Then, $\voffen \disjunion \leeremenge$ and $\compl\vkomp \disjunion \raumb$
are disjoint open neighbourhoods of $\elrauma$ and $\elraumb$, respectively.
\erueck
\bhin
Lemma \ref{lem:rel-top} implies that $\rauma$ and $\raumb$ are Hausdorff. 
To prove local compactness, let $\elrauma \in \rauma$ be given. 
Note first, that $\elrauma$ and $\twistabb(\elrauma)$ cannot be separated
using a type-1 neighbourhood of $\twistabb(\elrauma)$, 
as this point is never contained in such a set.
Note second, that both points can neither be separated by 
type-$\typzd$ sets. In fact, 
\bgl
\text{$\elrauma \in \bigl(\compl \vkomp_1 \cap \twistabb^{-1}(W_1)\bigr) \sqcup W_1$}
\quad
\text{and}
\quad
\text{$\twistabb(\elrauma) \in \bigl(\compl \vkomp_2 \cap \twistabb^{-1}(W_2)\bigr) \sqcup W_2$}
\egl
implies
$\twistabb(\elrauma) \in W_2$, hence 
$\elrauma \in \twistabb^{-1}(W_2)$, hence $\twistabb^{-1}(W_1 \cap W_2)$ is non-empty; contradiction.
As, on the other hand, 
two distinct elements in any Hausdorff
space have to be separable by elements of any basis,
$\elrauma$ and $\twistabb(\elrauma)$ are now separated by 
a type-1 and a type-$\typzd$ set, i.e., by
\bgl
\text{$\voffen \disjunion \leeremenge$}
\quad
\text{and}
\quad
\text{$\bigl(\compl \vkomp \cap \twistabb^{-1}(W)\bigr) \sqcup W$} \,,
\egl
respectively. As then $\voffen$ and $\compl \vkomp \cap \twistabb^{-1}(W)$
are disjoint, we see that $\voffen \cap \twistabb^{-1}(W)$ is contained in $\vkomp$.
Now we are done, as the first set is an open neighbourhood of $\elrauma$ and
the latter one is compact.
\qed
\ehin
\epf

\bprop
$\stdtwisttop$
is compact iff $\raumb$ is compact.
\eprop

\bpf
As $\raumb$ is always closed in $\stdtwisttop$, the compactness of 
$\stdtwisttop$ implies that of $\raumb$. Let us now prove the 
other direction assuming $\raumb$ to be compact.
\bunum
\item
Let $\cover$ be an open cover of $\stdtwisttop$.
We may assume that $\cover$ is contained in the basis of the topology.

\item
As $\cover$ covers, in particular, the compact set
$\raumb \teilmenge \stdtwisttop $, there is a finite 
$\subcover \teilmenge \cover$ still covering $\raumb$.
We may assume that none of the sets in $\subcover$ 
is of first type, as their intersection with $\raumb$ is empty.
Hence the
elements of $\subcover$ 
are type-$\typzd$ sets
\bgl
\elsubcover_i 
 \breitrel =  \bigl(\compl K_i \cap \twistabb^{-1}(W_i)\bigr) \sqcup W_i \,.
\egl
As $\subcover$ covers $\raumb$, we have $\bigcup W_i = \raumb$.
Now, for $K := \bigcup K_i$, we have
\bglklein[1.5ex]
\rauma \cap \: \bigcup_i \elsubcover_i 
 & = & \bigcup_i\bigl(\compl K_i \cap \twistabb^{-1}(W_i)\bigr) 
 \breitrel \obermenge 
        \compl K \cap \twistabb^{-1}(\bigcup_i W_i) 
 \breitrel= \compl K \,.
\eglklein
This means that $\subcover$ covers $\compl K \disjunion \raumb$.
\item
As $K$ is compact, we may find some finite $\subcover' \teilmenge \cover$ covering
$K$. Now, $\subcover \cup \subcover'$ covers all of $\stdtwisttop$.
\qed
\eunum
\epf
Let us finally compare the twisted-sum topology on $\rauma \disjunion \raumb$
with the standard direct-sum topology thereon.
As any standard set is open in the direct-sum topology, we have

\blem
The twisted-sum topology is always contained in the direct-sum topology.
\elem
\noindent
Criteria for the equality of both sum topologies are summarized in

\bprop
\label{prop:twist-vs-disjsum-krit}
Consider the following statements:

\bnum3
\item
\label{item:twist=disj-sum}
The sum topologies on $\rauma \disjunion \raumb$ coincide.
\item
\label{item:Y-compact}
$\rauma$ is compact.
\item
\label{item:Y-locallycompact}
$\rauma$ is locally compact.
\item
\label{item:f(Y)-closed}
$\twistabb(\rauma)$ is closed in $\raumb$.
\item
\label{item:Z-ist-offen}
$\leeremenge \disjunion \raumb$ is open in $\stdtwisttop$.
\item
\label{item:Z-huebsch-ueberdeckt}
$\raumb$ can be covered by open $W_\alpha$, whose
preimages $\twistabb^{-1}(W_\alpha)$ are contained in compacta in $\rauma$.
\enum
These statements are correlated as follows:

\makeatletter
\newcommand{\txtinvimpliz}[2][]{\ext@arrow 0359\Leftarrowfill@{#1}{#2}}
\newcommand{\txtimpliz}[2][]{\ext@arrow 0359\Rightarrowfill@{#1}{#2}}
\newcommand{\txtaequ}[2][]{\ext@arrow 0359\Leftrightarrowfill@{#1}{#2}}
\newcommand{\hdtxt}{$\raumb$ Hausdorff}
\newcommand{\downtxt}{dumi}
\makeatother

\begin{equation}
\begin{minipage}[B]{0.7\textwidth}
\leerezeile
\leerezeile
\bgl
\text{\ref{item:twist=disj-sum}} 
 \txtaequ{\text{\phantom{\hdtxt}}}
 \text{\ref{item:Z-ist-offen}} 
 \txtaequ{\text{\phantom{\hdtxt}}}
 \text{%
 \rput[Bl](0,1.43){\text{\ref{item:Y-compact}}}%
 \rput[Bl](0,-1.43){\text{\ref{item:Y-locallycompact}}}%
 \ref{item:Z-huebsch-ueberdeckt}%
 \rput[r]{-90}(-0.05,0.35){$\txtimpliz{\phantom{\downtxt}}$}%
 \rput[l]{-90}(-0.05,-0.12){$\txtimpliz{\phantom{\downtxt}}$}%
 }
 \txtimpliz{\text{\hdtxt}}
 \text{\ref{item:f(Y)-closed}} \\
\egl
\leerezeile
\leerezeile
\leerezeile
\leerezeile
\end{minipage}\nonumber
\end{equation}


\eprop

\bpf
\bvl{}
\iitem[\ref{item:twist=disj-sum} $\impliz$ \ref{item:Z-ist-offen}]
Trivial.
\iitem[\ref{item:Z-ist-offen} $\impliz$ \ref{item:twist=disj-sum}]
Observe 
\bgl
V \disjunion W 
 & = & (V \disjunion \leeremenge) \cup \bigl[(\leeremenge \disjunion \raumb) \cap (\twistabb^{-1}(W) \disjunion W)\bigr]
\egl
for any $V \teilmenge \rauma$ and $W \teilmenge \raumb$.
\iitem[\ref{item:Z-ist-offen} $\impliz$ \ref{item:Z-huebsch-ueberdeckt}]
As $\leeremenge \disjunion \raumb \teilmenge \stdtwisttop$
contains only the trivial type-1 set, it is a union of some type-$\typzd$ sets
\bgl
\bigl(\compl\vkomp_\alpha \cap \twistabb^{-1}(W_\alpha)\bigr) \disjunion W_\alpha \,.
\egl
By construction, each $\compl\vkomp_\alpha \cap \twistabb^{-1}(W_\alpha)$ is empty, i.e., 
$\twistabb^{-1}(W_\alpha) \teilmenge \vkomp_\alpha$.
As the $W_\alpha$ form an open cover of $\raumb$, we get the proof from continuity of $\twistabb$.
\iitem[\ref{item:Z-huebsch-ueberdeckt} $\impliz$ \ref{item:Z-ist-offen}]
This is clear from the proof of the reversed implication.
\iitem[\ref{item:Y-compact} $\impliz$ \ref{item:Z-huebsch-ueberdeckt}]
Trivial.
\iitem[\ref{item:Z-huebsch-ueberdeckt} $\impliz$ \ref{item:Y-locallycompact}]
Given such $W_\alpha$, there are compact
$\vkomp_\alpha$ containing the open $\twistabb^{-1}(W_\alpha)$.
Therefore, $\vkomp_\alpha$ is a compact neighbourhood for all the points 
in $\twistabb^{-1}(W_\alpha)$. As the latter sets form a cover
of $\rauma$, we get the claim.
\iitem[\ref{item:Z-huebsch-ueberdeckt} $\impliz$ \ref{item:f(Y)-closed}]
Assume that $\twistabb(\rauma)$ is not closed, i.e., there is some
$\elraumb \in \randrest{\twistabb(\rauma)}$.
Choose any open $W_\alpha =: W$ containing $\elraumb$. 
Because $W$ is open, there is%
\footnote{With $A := \twistabb(\rauma)$, we have $\elraumb \notin A \cap W$.
But, for each open $U$ containing $\elraumb$, 
we get an open $W \cap U$ containing $\elraumb$,
whence $\elraumb \in \quer A$ implies $\elraumb \in (A \cap W) \cap U$.
Thus, $\elraumb \in \quer{A \cap W}$, whence there is a net in $A \cap W$
converging to $\elraumb$.}
a net $(f(\elrauma_\lambda))$ in $\twistabb(\rauma) \cap W$ 
converging to $\elraumb$.
Since $\twistabb^{-1}(W)$ is contained in some compactum $\vkomp$,
there is a subnet of $(\elrauma_\lambda)$ converging to some 
$\elrauma \in \vkomp$. Consequently, a subnet of
$\twistabb(\elrauma_\lambda)$ converges to $\twistabb(\elrauma) \in \twistabb(\rauma)$.
As $\raumb$ is Hausdorff, we get $\twistabb(\elrauma) = \elraumb$,
hence a contradiction.
\qed
\evl
\epf

\noindent
None of the implications  
in Proposition \ref{prop:twist-vs-disjsum-krit} above
can be reversed, in general, nor can the Haus\-dorff property be removed there.
Let us explain this by means of several examples.

\bex
\label{ex:Z-grobe-top}
Let $\raumb$ carry the coarsest topology. 

Then, $\{\raumb\}$ is the only open cover for $\raumb$.
Proposition \ref{prop:twist-vs-disjsum-krit} now implies 
that both sum topologies coincide 
iff
$\rauma = \twistabb^{-1}(\raumb)$ is compact.
In particular, there are locally compact
$\rauma$ for which the two sum topologies differ, i.e.,
\ref{item:Y-locallycompact} $\impliz$ \ref{item:Z-huebsch-ueberdeckt}
is not given. 

Moreover,
taking any constant $\twistabb$,
we see that the image 
$\twistabb(\rauma)$ is closed iff $\raumb$ consists of a 
single point only or, equivalently, $\raumb$ is Hausdorff. 
Thus, for any compact $\rauma$
(implying the desired equality of
the twisted and the direct-sum topologies)
and any non-Hausdorff $\raumb$, 
\ref{item:Z-huebsch-ueberdeckt} $\impliz$ \ref{item:f(Y)-closed}
is not given.
In other words, the Hausdorff property is necessary.
\eex

\bex
Let $\rauma = \raumb$ be Hausdorff spaces and let $\twistabb$ be the identity.

We are going to show that the two sum
topologies coincide iff 
$\rauma$ is locally compact. The ``only if''-part is already covered by the 
proposition above. To show the ``if''-part, observe that,
for each $\elrauma \in \rauma$, we find open $W_\elrauma$ and closed compact $\vkomp_\elrauma$
with $\elrauma \in W_\elrauma \teilmenge \vkomp_\elrauma$. 
The claim follows from $\twistabb^{-1}(W_\elrauma) \ident W_\elrauma$
and Proposition \ref{prop:twist-vs-disjsum-krit}. 
Altogether, 
\ref{item:Z-huebsch-ueberdeckt} $\impliz$ \ref{item:Y-compact} is not given.

On the other hand, $\twistabb(\rauma) \ident \raumb$ is always closed in $\raumb$,
but $\raumb$ need not be locally compact, whence there is no need for the
twisted and the direct sum topologies to coincide. In other words, 
\ref{item:f(Y)-closed} $\impliz$ \ref{item:Z-huebsch-ueberdeckt} 
is not given. This remains true even if the Hausdorff assumption is dropped.

\eex

\bex
Let $\rauma = \raumb \kreuz \raumb$
and 
$\twistabb : \raumb \kreuz \raumb \nach \raumb$ be the projection
to the first component.

Of course, any $W \teilmenge \raumb$ has $W \kreuz \raumb$ as preimage,
which surely is contained in some compactum iff $\raumb$ is compact (or $W$ is empty). Therefore,
also
\ref{item:f(Y)-closed} and \ref{item:Y-locallycompact}
together
(even if $\raumb$ was Hausdorff)
do not suffice to imply 
\ref{item:Z-huebsch-ueberdeckt} 
\eex

\newpage

\section{Preliminaries on $C^\ast$-algebras}
\label{sect:c-ast-algs}

Before going to the main statements,
let us summarize the relevant prerequisites from $C^\ast$-algebras.
Note that we assume any ideal to be closed and two-sided.

\subsection{Closedness of Subalgebra Sums}

For completeness, let us start with the well-known \cite{Murphy} 
\bprop
\label{prop:ideal+subalg=Cstern}
Let $\clg$ be a $C^\ast$-algebra with ideal $\ideal$ and $C^\ast$-subalgebra $\alg$.

Then $\ideal + \alg$ is closed, hence a $C^\ast$-subalgebra of $\clg$.
\eprop
\bpf
As $\ideal$ is an ideal, $\clg/\ideal$ is a $C^\ast$-algebra and
the canonical projection $\pi : \clg \nach \clg/\ideal$ is a $\ast$-homomorphism.
It restricts to a $\ast$-homomorphism
$\pi_{\alg} : \alg \nach \clg/\ideal$. 
Consequently \cite{Murphy}, 
the range of $\pi_{\alg}$ is closed, hence $\pi^{-1}(\pi_{\alg}(\alg)) = \ideal + \alg$ 
as well.
\qed
\epf

\noindent
Sometimes, it might not be clear a priori whether $\ideal$ is indeed an ideal
in $\clg$ -- or even worse what $\clg$ really is.
This, however, does not destroy the closedness of $\ideal + \alg$
as long as at least the relation between $\ideal$ and $\alg$ resembles
the ideal property:

\bcorr
\label{corr:quasiideal+subalg=Cstern}
Let $\mlg$ and $\alg$ be $C^\ast$-subalgebras of some $C^\ast$-algebra $\clg$ with 
$\alg \mlg \teilmenge \mlg$. 

Then $\allgsumalg$ is a $C^\ast$-subalgebra of $\clg$, containing $\mlg$ as an ideal.
\ecorr
\bpf
Obviously, $\blg := \allgsumalg$ is a $\ast$-subalgebra of $\clg$ with
$\blg \mlg \teilmenge \mlg$.
Consequently, its closure $\quer \blg$ is a $C^\ast$-subalgebra of $\clg$.
Now, $\mlg$ and $\alg$ are also $C^\ast$-subalgebras of $\quer\blg$.
Even more, $\mlg$ is an ideal there. In fact, given $\elclg \in \quer\blg$,
there are $\elblg_i$ in $\blg$ converging to $\elclg$. For any
$\elmlg \in \mlg$, we have now $\elblg_i \elmlg \in \mlg$,
whence $\elclg \elmlg = \lim_i \elblg_i \elmlg \in \mlg$.
By Proposition \ref{prop:ideal+subalg=Cstern}, $\allgsumalg$ is
closed in $\quer\blg$, hence equal to $\quer\blg$ by denseness.
\qed 
\epf
In other words, replacing $\clg$ above by $\allgsumalg$ returns ourselves to the
previous situation.
Therefore, we will restrict ourselves to the case that $\mlg$ is an ideal.

%

\subsection{Gelfand Transform}

Let $\alg$ be an abelian $C^\ast$-algebra with spectrum $\malg$.
Recall
\cite{Murphy}

\bdf
\bnum2
\item
The \df{Gelfand transform} $\gelf \elalg$
of any $\elalg \in \alg$ is given by
\bgl
\fktdefabgesetzt{\gelf\elalg}{\malg}{\C\,.}{\charact}{\charact(\elalg)}
\egl
\item
The topology on $\malg$ is the initial topology 
induced by all the Gel\-fand transforms.
More precisely, it is generated by all 
the sets $\gelf\elalg^{-1}(\offenCalg)$
with $\elalg \in \alg$ and open $\offenCalg \teilmenge \C$.
\enum
\edf

\bprop
\label{prop:inittop-erzeugt-von-teilmenge}
Let $\blg \teilmenge \alg$ be any subset that generates $\alg$ as a
$C^\ast$-algebra.

Then the topology on $\malg$ is 
already
induced by $\{\gelf \elalg \mid \elalg \in \blg\}$.
\eprop
\noindent
Note that the Gelfand transform of $\elalg \in \blg$ is taken w.r.t.\ $\alg$,
not w.r.t.\ $\blg$.

For the proof of the proposition above, recall that the initial topology on
some topological space $\dummyraumb$, induced by some functions $\dummyfktb_\alpha$
thereon,
is characterized by the following condition: For any map 
$\dummyfkt : \dummyraum \nach \dummyraumb$ on
any topological space $\dummyraum$, we have
\bgl
\text{$\dummyfkt$ continuous} 
& \aequ & 
\text{$\dummyfktb_\alpha \circ \dummyfkt$ continuous for all $\alpha$.}
\egl
\noindent
Now, the proposition above
is an immediate consequence of
\blem
Let $\dummyraum$ be a topological space and $\dummyfkt : \dummyraum \nach \malg$.

Then 
$\gelf\elalg \circ \dummyfkt$ is continuous for all $\elalg \in \alg$
iff 
$\gelf\elblg \circ \dummyfkt$ is continuous for all $\elblg \in \blg$. 
\elem
\bpf
By continuity of the algebra operations (sum, product, conjugation),
we may assume $\alg = \quer\blg$. Moreover, we only have to prove the ``if''-part.
For this, write $\elalg \in \alg$ as $\elalg = \lim \elblg_k$ with 
$\elblg_k \in \blg$. Then 
\bgl
\supnorm{\gelf\elblg_k \circ \dummyfkt - \gelf\elalg \circ \dummyfkt}
 \breitrel\leq \supnorm{\gelf\elblg_k - \gelf\elalg}
 \breitrel= \norm{\elblg_k - \elalg}
 \breitrel\gegen 0
\egl
by linearity and isometry of the Gelfand transform. Consequently,
$\gelf\elalg \circ \dummyfkt$ is continuous.

\qed
\epf

\subsection{Natural Mapping}

Let $X$ be some locally compact Hausdorff space and let
$\boundfkt(X)$ be the $C^\ast$-algebra of bounded continuous functions on it.
Moreover, let $\alg$ be a unital $C^\ast$-subalgebra of $\boundfkt(X)$. 
Recall from \cite{e8,paper39} the following definition and proposition.
\bdf
The \df{natural mapping}
$\iota:\set \nach \malg$ is given by
\bgl
\fktdefabgesetzt{\iota(\elset)}{\alg}{\C\,.}{\elalg}{\elalg(\elset)} 
\egl
\edf

\bprop
\label{prop:rendall}
\bunum
\item 
$\iota$ is well defined and continuous. 
\item 
$\iota$ is injective iff $\alg$ separates the points in $\set$. 
\item 
$\iota(\set)$ is dense in $\malg$. 
\item 
$\gelf \elalg \circ \iota = \elalg$ on $\set$ for all
$\elalg \in \alg$.
\eunum 
\eprop
%
%


\section{Topology of $\malgab$}
\label{sect:sumspec_top}

From now on, we will use the following

\bnotat
\bnotatlist
\citem[X] some nonempty locally compact Hausdorff space
\citem[\alga] some ideal in $C_0(X)$
\citem[Y] the open subset of $X$ with 
\bgl
\alga & = & \cnulleinschr Y X \breitrel{:=} \bigl\{\elalga \in C_0(X) \mid \text{$\elalga \ident 0$ on $\compl Y$}\bigr\}
\egl%
\citem[\algb] some unital $C^\ast$-subalgebra of $\boundfkt(X)$
with
\bgl
\alga \cap \algb & = & \NULL 
\egl
\citem[\algab] the direct vector space sum 
\bgl
\algab & := & \algablang 
\egl
\citem[\twist] the restriction of the natural mapping $\iotaalgb : X \nach \malgb$
to $Y$:%
\bgl
\twist & := & \iotaalgb\einschr Y : Y \nach \malgb 
\egl

\enotatlist
\enotat

\brem
\bnum2
\item
We assume the direct vector space sum of subspaces of a third vector space
to be contained in this vector space again.
\item
Note that $Y$ is well defined as any ideal in $C_0(X)$ is of the form
$\cnulleinschr Y X$.~\cite{ConwayFktA} 
Moreover, as $\cnulleinschr Y X$ is naturally isomorphic to $C_0(Y)$, 
we will identify its spectrum with $Y$.
\enum
\erem
\enlargethispage{0.5\baselineskip}

\blem
\bnum2
\item
$\algab$ is a unital abelian $C^\ast$-subalgebra of $\boundfkt(X)$.
\item
$\alga$ is an ideal in $\algab$.
\enum
\elem
\bpf
\bnum2
\item
$\alga \ident \cnulleinschr Y X$ is an ideal in $\boundfkt(X)$, whence 
the statement follows from
$\algab \ident \alga + \algb$ and 
Proposition \ref{prop:ideal+subalg=Cstern}.
\item
This is trivial by the preceding argument.
\qed
\enum
\epf

\noindent
We are now stating our main result.

\bthm
\label{thm:spec_sum}
We have
\bgl
\malgab & \iso & \twisttop\Y\malgb\twist \,.
\egl
\ethm

\bpf
Let us define the mapping 
\bgl
\charabb \breitrel: \twisttop\Y\malgb\twist \nach \malgab
\egl
for $y \in \Y$ and $\charaktalgb \in \malgb$ as follows:
\bgl{}
[\charabb(y)] (\elalga + \elalgb) & := & (\elalga + \elalgb) (y) \\{}
[\charabb(\charaktalgb)] (\elalga + \elalgb) & := & \charaktalgb(\elalgb) \,.
\egl
We may assume $\Y$ to be non-empty.
\bunum
\item
\df{Well-definedness}

Of course, $\charabb(y) \in \malgab$. For the other part observe that 
$\charabb(\charaktalgb)$ is nonzero and
multiplicative by 
\bgl{}
[\charabb(\charaktalgb)]\bigl((\elalga + \elalgb)(\elalga' + \elalgb')\bigr)
  & = & [\charabb(\charaktalgb)]\bigl((\elalga \elalga' + \elalga \elalgb' + \elalgb \elalga') + \elalgb \elalgb'\bigr) \\
  & = & \charaktalgb(\elalgb \elalgb') \\
  & = & \charaktalgb(\elalgb) \: \charaktalgb(\elalgb') \\
  & = & [\charabb(\charaktalgb)](\elalga + \elalgb) \: [\charabb(\charaktalgb)](\elalga' + \elalgb').
\egl
\item
\df{Surjectivity}

Let $\charakt : \algab \nach \C$ be a character on $\algab$. Then
there are two cases:
\bunum
\item
If $\charact \einschr{\alga} = 0$, 
then, obviously, $\charact\einschr{\algb}$ is a character on $\algb$,
with $\charabb(\charact\einschr{\algb}) = \charact$.
\item
If $\charact \einschr{\alga} \neq 0$, then it is a character on $\alga$,
whence,
by Gelfand-Naimark theory,
there is some $\ely \in \Y$ with $\charact(\elalga) = \elalga(\ely)$ for all $\elalga \in \alga$. 
Given $\elalgb \in \algb$, we have for some $\elalga$ with $\charact(\elalga) \neq 0$
\bgl
\elalga(\ely) \elalgb(\ely) \breitrel= (\elalga \elalgb) (\ely) \breitrel= \charact(\elalga \elalgb) 
\breitrel= \charact(\elalga) \charact(\elalgb) \breitrel= \elalga(\ely) \charact(\elalgb),
\egl
whence $\charact(\elalgb) = \elalgb(\ely)$. Here, we used $\alga \algb \teilmenge \alga$. 
Thus, we have $\charact(\elalgab) = \elalgab(\ely)$
for all 
$\elalgab = \elalga + \elalgb \in \algab$, hence $\charact = \charabb(\ely)$.
\eunum
\item
\df{Injectivity}

There are three cases:
\bunum
\item
Let $\ely,\ely' \in \Y$ with $\ely \neq \ely'$. Taking some $\elalga \in \alga$ with
$\elalga(\ely) \neq \elalga(\ely')$, we get
\bgl\
[\charabb(\ely)](\elalga) 
  \breitrel= \elalga(\ely) 
  \breitrel\neq \elalga(\ely') 
  \breitrel= [\charabb(\ely')](\elalga)\,,
\egl
implying
$\charabb(\ely) \neq \charabb(\ely')$.
\item
Let $\ely \in \Y$ and $\charaktalgb \in \spec \algb$. Then, for any $\elalga \in \alga$ with
$\elalga(\ely) \neq 0$, we have
\bgl\
[\charabb(\ely)](\elalga) 
  \breitrel= \elalga(\ely) 
  \breitrel\neq 0 
  \breitrel= [\charabb(\charaktalgb)](\elalga) \,,
\egl
implying 
$\charabb(\ely) \neq \charabb(\charaktalgb)$.
\item
Let $\charaktalgb, \charaktalgb' \in \spec \algb$ with $\charaktalgb \neq \charaktalgb'$. Take
$\elalgb \in \algb$ with $\charaktalgb(\elalgb) \neq \charaktalgb'(\elalgb)$.
Then 
\bgl\
[\charabb(\charaktalgb)](\elalgb) 
  \breitrel= \charaktalgb(\elalgb) 
  \breitrel\neq \charaktalgb'(\elalgb) 
  \breitrel= [\charabb(\charaktalgb')](\elalgb) \,,
\egl
implying 
$\charabb(\charaktalgb) \neq \charabb(\charaktalgb')$.
\eunum
\item
\df{Continuity}

By Proposition \ref{prop:inittop-erzeugt-von-teilmenge}, the topology of $\malgab$ is generated by the sets 
$\widetilde \elalga^{-1} (\offenCalga)$ and $\widetilde \elalgb^{-1} (\offenCalgb)$ 
with $\elalga \in \alga$, $\elalgb \in \algb$ and open $\offenCalgboth \teilmenge \C$.
So, we only have to show that their preimages are open.

\bunum

\item
Let $\elalga \in \alga$.
Then 
\bgl{}
[\widetilde \elalga \circ \charabb] (\ely) 
 & = & [\charabb(\ely)](\elalga) 
 \breitrel= \elalga(\ely) 
\\{}
[\widetilde \elalga \circ \charabb] (\charaktalgb) 
 & = & [\charabb(\charaktalgb)](\elalga) 
 \breitrel= 0 
\egl
for $\ely \in \Y$ and $\charaktalgb \in \spec \algb$. Let now $\offenCalga \teilmenge \C$ be open.
\bunum
\item
If $0$ is not contained in $\offenCalga$, then 
\bgl\
\charabb^{-1} \bigl(\widetilde \elalga^{-1}(\offenCalga)\bigr) 
 \breitrel\ident (\widetilde \elalga \circ \charabb)^{-1}(\offenCalga)
 \breitrel= \elalga^{-1}(\offenCalga)
 \breitrel\ident \elalga^{-1}(\offenCalga) \disjunion \leeremenge
\egl
is a type-1 element, as $\elalga \in \alga \ident C_{0,\Y}(X)$,
whence $\elalga^{-1}(\offenCalga)$ is open in $\Y$.
\item
If $0$ is contained in $\offenCalga$, then 
\bgl\
\charabb^{-1} \bigl(\widetilde \elalga^{-1}(\offenCalga)\bigr) 
 \breitrel\ident (\widetilde \elalga \circ \charabb)^{-1}(\offenCalga)
 \breitrel= \elalga^{-1}(\offenCalga) \disjunion \malgb
\egl
This is a type-2 element
since
the complement of $\elalga^{-1}(\offenCalga)$ in $\Y$ is
compact, for $\elalga \in \alga$.

\eunum
\item
Let $\elalgb \in \algb$. Then, by Proposition \ref{prop:rendall},
\bgl{}
[\widetilde \elalgb \circ \charabb] (\ely) 
 & = & [\charabb(\ely)](\elalgb) 
 \breitrel= \elalgb(\ely) 
 \breitrel= \widetilde\elalgb(\iota_1(\ely))
 \breitrel\ident \widetilde\elalgb(\twist(\ely))
\\{}
[\widetilde \elalgb \circ \charabb] (\charaktalgb) 
 & = & [\charabb(\charaktalgb)](\elalgb) 
 \breitrel= \charaktalgb (\elalgb)
 \breitrel= \widetilde \elalgb (\charaktalgb)
\egl
for $\ely \in \Y$ and $\charaktalgb \in \spec \algb$. 
This means that
\bgl
\charabb^{-1} \bigl(\widetilde \elalgb^{-1} (\offenCalgb)\bigr) 
 \breitrel\ident [\widetilde \elalgb \circ \charabb]^{-1} (\offenCalgb)
 \breitrel= \twist^{-1}(\gelf\elalgb^{-1}(\offenCalgb)) \disjunion \gelf\elalgb^{-1}(\offenCalgb)
\egl
is a type-3 element as 
$\gelf\elalgb^{-1}(\offenCalgb)$ is open in $\malgb$.
\eunum
\item
\df{Homeomorphy}

As $\charabb$ is a continuous bijection from a compact to a Hausdorff space,
it is even a homeomorphism.
\qed
\eunum
\epf


\section{Embedding and Denseness}
\label{sect:denseness}

From the proof of Theorem \ref{thm:spec_sum}, we get immediately
\bcorr
\label{corr:spec-sum_nat-mapp}
The natural mapping 
$\iotaalgab : X \nach \malgab$ is given by
\bgl
\iotaalgab & = & \begin{cases}
                \charabb & \text{on $\Y$} \\
                \charabb \circ \iotaalgb & \text{on $X \setminus \Y$} \,.
               \end{cases} 
\egl
Here, $\iotaalgb : X \nach \malgb$ is the natural mapping w.r.t.\ $\algb$.
\ecorr
Ignoring the homeomorphism $\charabb : Y \disjunion \malgb \nach \malgab$, the
natural mapping for $\algab$ equals the identity on $Y$ and
the natural mapping $\iotaalgb$ for $\algb$ on $Y \setminus X$.
In particular, for $Y = X$, 
the natural mapping is simply the identity on $X$.

\blem
$X$ is densely and continuously embedded into $\malgab$, 
provided $\algb$ separates the points in $X \setminus Y$.
\elem
\bpf
First of all, $\alga \ident \cnulleinschr Y X$ separates any point in $Y$ from
any other point in $X$.
Next, by assumption, $\algb$ separates any two points in $X \setminus Y$.
Altogether, $\alg$ separates any two points in $X$.
The statement follows from Proposition \ref{prop:rendall}.
\qed
\epf

\bcorr
Let $X \setminus \Y$ contain at most one point.

Then 
$X$ is densely and continuously embedded into $\malgab = \twisttop\Y\malgb\twist$.
\ecorr
Note that, for $\Y = X$,
the space $X$ might be embedded into $\malgab$ in two different ways: 
\bunum
\item
Firstly, one uses the natural embedding $\iota$ of $X$ into $X \disjunion \malgb$,
whose image is again $X$ seen as a subset of $X \disjunion \malgb$. That is the way we 
went above.
\item
Secondly, assuming that $\algb$ separates the points in $X$,
one uses the natural embedding $\iotaalgb : X \nach \malgb$ and then
embeds $\malgb$ into $X \disjunion \malgb$.
However, now the resulting embedding of $X$ into $X \disjunion \malgb$ 
is not dense anymore, unless $X$ is empty. This is clear, as
$\malgb$ is closed in $\malgab$.

\eunum

\noindent
A similar behaviour can be observed if $\Y$ equals $X$ minus some point. Then
the natural mapping $\iotaalgab$ is given by the identity on $\Y$, but the ``missing'' point
is taken from $\malgb$. In other words, $\malgb$ is attached to $\Y$ ``filling'' the gap.


\section{Examples}
\label{sect:examples}

\bex
\label{ex:1pt-comp-explizit}
\df{One-point compactification}

Let 
$\alga = C_0(X)$ and consider $\algb := \C\,\EINS$ as a subset of $\boundfkt(X)$.
As $\malgb$ consists of a single point, say $\infty$, only, we have exactly two open sets:
$\leeremenge$ and $\{\infty\}$. Moreover, the twisting map $\twist : X \nach \malgb$ is
trivial. Consequently, the only type-3 sets are
$\leeremenge \disjunion \leeremenge$ and $X \disjunion \{\infty\}$.
This means that the topology of
$\spec(C_0(X) \dirvsum \C\,\EINS)$ 
is generated just by
the open sets in $X$ and 
by the complements of closed compact sets in $X$ united with $\spec \algb$.
This is indeed nothing but the topology of the one-point compactification $X^\ast$ of $X$.
Of course, $X$ is dense in $X^\ast$.
Moreover, Lemma \ref{lem:hausdorff-krit} generalizes the well-known
fact \cite{Kelley}, 
that $X^\ast$ is Hausdorff iff $X$ is locally compact Hausdorff.
Finally, Example \ref{ex:Z-grobe-top}
comprises the fact \cite{Kelley} 
that $\infty$ is an isolated point of $X^\ast$
iff $X$ is compact.

\eex

\bex
\label{ex:aap}
\df{Asymptotically almost periodic functions}

Let $X$ be a non-compact, but locally compact abelian group, and let $\alga$
be full $C_0(X)$.
If $\algb$ is the set of almost periodic functions on $X$, then $\algab$ is the 
set of asymptotically almost periodic functions. (See \cite{GrigTonev}
for $X = \R$.)
Its spectrum is given by the twisted sum 
\bgl
\twisttop X{X_\Bohr}\iota \,,
\egl
where 
$X_\Bohr \ident \malgb$ is the Bohr compactification \cite{RudinFour} of $X$,
and $\iota$ is the canonical embedding of $X$ into $X_\Bohr$.
Open sets are, in particular, the open sets in $X$ and
the type-3 sets $f^{-1}(U) \disjunion \gelf f^{-1}(U)$ with open $U \teilmenge \C$ 
and with $f$ 
running over the almost periodic functions.
As $\iota(X)$ is dense in the compactum 
$X_\Bohr$ (being, e.g., a consequence of Proposition \ref{prop:rendall}),
but strictly smaller, it cannot be a closed subset. 
Now, Proposition \ref{prop:twist-vs-disjsum-krit} implies that $\twisttop X{X_\Bohr}\iota$ is not the
direct-sum topology on $X \disjunion X_\Bohr$.
\eex


\section{Acknowledgements}
The author thanks Maximilian Hanusch
for numerous discussions and many helpful comments on 
a draft of the present article.
Moreover, the author gratefully acknowledges 
discussions with
Thomas Tonev concerning asymptotically almost periodic functions.
The author has been supported by the Emmy-Noether-Programm of
the Deutsche Forschungsgemeinschaft under grant FL~622/1-1.


\anhangengl


\section{Borel algebra of $\stdtwisttop$}
\label{app:measure}

Reusing the notation from Section \ref{sect:twisted-top}, we will now describe
the Borel algebra that corresponds to the twisted sum, as well as 
the (regular) measures on it. For the particular case of asymptotically
almost periodic functions see also \cite{max-diss}. 

\bprop
\label{prop:dirsum-borelalg}
We have 
\bgl
\borelalg(\stdtwisttop) & = & \borelalg(\rauma) \dirsum \borelalg(\raumb) \,.
\egl
\eprop
\noindent
Here, $\borelalg(X)$ denotes the Borel $\sigma$-algebra of a topological space $X$.
The proposition above shows, in particular, 
that the direct sum and the twisted sum topologies yield
one and the same Borel algebra.
For brevity, we refrain
from indicating $\twistabb$ in $\stdtwisttop$ in the following.
\simplifytwist

\bpf
\bunum
\item
As $\rauma$ is open, any open set in $\rauma$ is open in $\stdtwisttop$, whence 
\bgl
\borelalg(\rauma) & \teilmenge & \borelalg(\stdtwisttop) \,.
\egl
As $\raumb$ is closed, any closed set in $\raumb$ is
closed in $\stdtwisttop$, whence 
\bgl
\borelalg(\raumb) & \teilmenge & \borelalg(\stdtwisttop) \,.
\egl
\item
Obviously,
$\borelalg(\rauma) \dirsum \borelalg(\raumb)$ is
a Borel algebra and contains any standard open set in $\stdtwisttop$.
As these sets generate the twisted topology, hence the respective Borel algebra,
we have 
\bgl
\borelalg(\rauma) \dirsum \borelalg(\raumb) & \obermenge & \borelalg(\stdtwisttop) \,.\phantom{{}\dirsum \borelalg(\raumb)}
\egl
\qed
\eunum
\epf

\newcommand{\maszsum}{\masza \dirsum \maszb}
\newcommand{\maszab}{\mu}
\newcommand{\masza}{\mu_1}
\newcommand{\maszb}{\mu_2}
\newcommand{\setrauma}{Y}
\newcommand{\setraumb}{Z}
\newcommand{\setraumab}{\setrauma \disjunion \setraumb}

\noindent
Recall that the direct sum $\maszsum$ of Borel
measures $\masza$ on $\rauma$ and $\maszb$ on $\raumb$ is given by
\bgl\
[\maszsum](\setrauma\disjunion\setraumb) 
  & := & \masza(\setrauma) + \maszb(\setraumb)
\qquad
\text{for $\setrauma \in \borelalg(\rauma)$ and $\setraumb \in \borelalg(\raumb)$\,.}
\egl
\noindent
Of course, $\maszsum$ is a Borel measure on 
$\stdtwisttop$ by Proposition \ref{prop:dirsum-borelalg}.
Even more, we have

\bcorr
The finite Borel measures on $\stdtwisttop$ are precisely the direct sums of 
finite Borel measures on $\rauma$ and on $\raumb$. The respective statements 
hold for measures that are additionally inner or outer regular.
\ecorr

\bpf
\bunum
\item
Obviously, any direct sum of finite Borel measures is finite Borel. The other way round,
given a finite Borel measure $\maszab$ on $\stdtwisttop$, then 
$\masza := \maszab \einschr{\borelalg(\rauma)}$ defines a finite Borel measure on
$\rauma$; similarly, $\maszb$ is constructed on $\raumb$.
They fulfill $\maszab = \masza \dirsum \maszb$.
\item
As any finite inner regular Borel measure is outer regular \cite{elstrodt}, 
we may restrict ourselves to the
inner regular case. First, on the one hand, let $\masza$ and $\maszb$ be inner regular and let
$\setrauma \disjunion \setraumb \in \borelalg(\stdtwisttop)$. Then, for $\maszab := \maszsum$, 
\bgl
\masza(\setrauma) + \maszb(\setraumb)
  & \ident & \maszab(\setrauma\disjunion\setraumb) \\
  & \geq & \sup\bigl\{\maszab(K) \mid \text{$K \teilmenge \setrauma\disjunion\setraumb$ compact} \bigr\} \\
  & \geq & \sup\bigl\{\maszab(K_1 \disjunion K_2)
               \mid \text{$K_1 \teilmenge \setrauma$, $K_2 \teilmenge \setraumb$ compact} \bigr\} \\
  & = & \sup\bigl\{\masza(K_1)
               \mid \text{$K_1 \teilmenge \setrauma$ compact} \bigr\} \\&&\hspace*{\fill}{}+
        \sup\bigl\{\maszb(K_2)
               \mid \text{$K_2 \teilmenge \setraumb$ compact} \bigr\} \\
  & = & \masza(\setrauma) + \maszb(\setraumb)\,,
\egl
whence $\maszsum$ is inner regular.
Here, the first inequality comes from monotonicity and the second one as any two compact
sets have compact union. Now, on the other hand,
let $\maszsum$ be inner regular. Then, however, both $\masza$ and $\maszb$ have to be
regular as well, since both spaces $\rauma$ and $\raumb$ carry the relative topology from 
$\stdtwisttop$, whence a subset of, say, $\rauma$ is compact iff if it is so in $\stdtwisttop$.
\qed
\eunum
\epf
\noindent
In particular, any measure on one of the subspaces $\rauma$ and $\raumb$ is
a measure on the twisted sum with the respective additional properties. 
So, in the case of asymptotically almost periodic (AAP)
functions, the Haar measure on $\R_\Bohr$ is also a measure on the AAP spectrum.
Moreover, according to Hanusch \cite{max-diss}, it is the only normalized regular Borel measure 
that is invariant w.r.t.\ the induced $\R$-action on $\twisttop\R{\R_\Bohr}\iota$.


\end{document}